\documentclass[12pt]{article}
\usepackage{amsmath}
\usepackage{amssymb}
\usepackage{amscd}
\newtheorem{theoreme}{Th\'eor\`eme}[section]
\newtheorem{definition}[theoreme]{D\'efinition}
\newtheorem{proposition}[theoreme]{Proposition}
\newtheorem{corollaire}[theoreme]{Corollaire}
\newtheorem{lemme}[theoreme]{Lemme}

\newtheorem{remarque}[theoreme]{Remarque}
\newenvironment{preuve}{\begin{trivlist} \item[]{\it Preuve---}}
{\par\hfill $\square$\end{trivlist}}
\renewcommand{\P}{\mathbb{P}}
\renewcommand{\L}{{\cal L}}
\newcommand{\C}{\mathbb{C}}
\newcommand{\R}{\mathbb{R}}
\newcommand{\N}{\mathbb{N}}
\newcommand{\Z}{\mathbb{Z}}
\newcommand{\T}{{\rm T}}

\renewcommand{\H}{{\cal H}}
\newcommand{\G}{{\cal G}}
\newcommand{\A}{{\cal A}}
\newcommand{\B}{{\cal B}}

\newcommand{\F}{{\cal F}}
\newcommand{\D}{{\cal D}}

\newcommand{\mult}{{\rm mult}}
\newcommand{\const}{{\rm const}}
\newcommand{\voir}{{\it voir }}
\newcommand{\Id}{{\rm Id}}
\newcommand{\ddc}{{\rm dd}^{\rm c}}
\renewcommand{\Re}{{\rm Re}}
\renewcommand{\Im}{{\rm Im}}
\renewcommand{\O}{{\rm O}}

\newcommand{\ahead}{\par \hspace{0.75cm}}
\title{Sur les applications de Latt\`es de $\P^k$
\footnote{{\bf Classification math\'ematique:} 30D05, 58F23.
\ \ \ \ \ \ \ \ \ \ \ \ \ \ \ \ \ \ \ \ \ \ \ \ \ \ \
\break
{\bf Mots cl\'es:} application de Latt\`es, ensemble de Julia,
courant de Green.}
}
\author{Tien-Cuong Dinh}
\begin{document}
\maketitle
\begin{abstract} Let $f$ be a polynomial endomorphism of 
degree $d\geq 2$ of $\C^k$ ($k\geq 2$)
which extends to a holomorphic endomorphism of $\P^k$. Assume that
the maximal order Julia set of $f$ is laminated by real
hypersurfaces in some open set.
We show that $f$ is homogenous and is a polynomial lift
of a Latt\`es endomorphism of $\P^{k-1}$.
\end{abstract}
{\it \underline{Notations}:}
\begin{itemize}
\item $f$ un endomorphisme de degr\'e $d\geq 2$ v\'erifiant
l'hypoth\`ese du th\'eor\`eme 1.3, $g$ un it\'et\'e de degr\'e
$d_g$ de $f$, $b$ un point fixe
r\'epulsif de $g$.
\item $G$ la fonction de Green, $\mu$ la mesure d'\'equilibre,
$J$ l'ensemble de Julia d'ordre maximal de $f$.
\item $\varphi$ une application de Poincar\'e de $g$ en $b$,
$\Lambda:=\varphi^{-1}\circ g\circ \varphi$,
$G^*:=G\circ\varphi$, $\mu^*:=\varphi^*(\mu)$,
$J^*:=\varphi^{-1}(J)$.
\item $z=(z_1,...,z_k)\in\C^k$, 
$z':=(z_1,...,z_{k-1})$ et $\|z'\|^2:=|z_1|^2+\cdots+|z_{k-1}|^2$.
\item Si $\G$ est un groupe d'automorphismes de $\C^k$ et $V$
un sous-ensemble de $\C^k$, on note $\G|_V:=\{\tau|_V \mbox{ avec }
\tau\in\G \mbox{ v\'erifiant } \tau(V)=\tau^{-1}(V)=V\}$.
\end{itemize}
\section{Introduction}
\ahead
Dans \cite{Lattes}, Latt\`es a construit des fractions
rationnelles semi-conjugu\'ees \`a des isog\'enies dilatantes
d'un tore complexe. L'ensemble de Julia d'une telle fraction
remplit la sph\`ere de Riemann et sa mesure d'\'equilibre est lisse
en dehors d'un ensemble fini.
\par
En dimension quelconque, Berteloot et Loeb ont caract\'eris\'e
les applications de Latt\`es de $\P^k$
par le fait que leurs courants de
Green sont lisses et strictement positifs dans un certain ouvert
\cite{BertelootLoeb2}.
Notons ici qu'en dimension 1, les notions de mesure
d'\'equilibre et de courant de Green sont co\"{\i}ncides.
Rappelons la d\'efinition d'applications de Latt\`es donn\'ee
par Berteloot et Loeb:
\begin{definition} \rm
Un endomorphisme holomorphe $f$ de degr\'e $d\geq 2$
de $\P^k$ est
appel\'e {\it application de Latt\`es}, s'il existe un groupe
d'isom\'etries complexes discret, co-compact $\A$ de $\C^k$, une
application affine $\Lambda_f$ de partie lin\'eaire $\sqrt{d}\tau$
($\tau$ unitaire) ainsi qu'un rev\^etement ramifi\'e
$\varphi:\C^k\longrightarrow \P^k$ tels que
$f\circ\varphi=\varphi\circ \Lambda_f$ et $\A$ agisse
transitivement sur les fibres de $\varphi$.
\end{definition}
\par
Dans \cite{DinhSibony}, en \'etudiant les endomorphismes
permutables de $\P^k$, nous avons introduit et \'etudi\'e une classe
d'endomorphismes dont la fonction de
Green, la mesure d'\'equilibre et l'ensemble de Julia sont
lamin\'es. Cette classe contient \'egalement les relev\'es
polynomiaux des applications de
Latt\`es. Les relations \'etroites
entre les fractions de Latt\`es et l'\'etude des fonctions
permutables ont \'et\'e \'eclair\'es par Eremenko dans un
travail ant\'erieur \cite{Eremenko} (\voir \'egalement
\cite{Fatou,Julia,Ritt}).
\par
Dans ce pr\'esent article, nous nous int\'eressons
\`a des endomorphismes polynomiaux dont l'ensemble de Julia
d'ordre maximal est, dans un certain ouvert,
lamin\'e par des hypersurfaces r\'eelles.
Nous allons prouver que ces endomorphismes sont homog\`enes et
leurs restrictions \`a l'hyperplan \`a l'infini sont des
applications de Latt\`es, c'est-\`a-dire qu'ils sont des
relev\'es polynomiaux des applications de Latt\`es.
Dans un cas particulier (pour les applications polynomiales
homog\`enes de $\C^{k+1}$
dont l'ensemble de Julia d'ordre maximal est, dans un
certain ouvert, une hypersurface r\'eelle strictement
pseudoconvexe) on obtient
la caract\'erisation des applications de Latt\`es
mentionn\'ee ci-dessus. Par m\^eme m\'ethode,
nous pouvons sans doute \'etendre notre
\'etude au cas o\`u l'ensemble de Julia d'ordre maximal est, dans
un certain ouvert, une
vari\'et\'e r\'eelle de petite codimension.
\par
Pour faciliter la lecture, nous allons rappeler quelques notions
fondamentales de la th\'eorie des syst\`emes dynamiques
holomorphes en
plusieurs variables et r\'esumer certaines propri\'et\'es utiles
pour cet article. Le lecteur trouve des expos\'es d\'etaill\'es
dans \cite{Sibony,Fornaess,BedfordSmillie}.
\par
Notons $\pi:\C^{k+1}\setminus \{0\} \longrightarrow \P^k$ la
projection canonique.
Soit $f$ un endomorphisme holomorphe de degr\'e $d\geq 2$
de $\P^k$. On note $f^s$ l'it\'er\'e $s$-i\`eme de $f$.
{\it Un relev\'e}
de $f$ est un endomorphisme polynomial homog\`ene de $\C^{k+1}$
v\'erifiant $\pi\circ F=f\circ\pi$. On appelle {\it fonction de
Green} de $f$ la fonction $G_f$ qui est d\'efinie par
$$G_f(z):=\lim_{s\rightarrow\infty}\frac{1}{d^s}\log\|F^s(z)\|.$$
C'est une fonction continue plurisousharmonique de
$\C^{k+1}\setminus\{0\}$ v\'erifiant $G_f\circ F=dG_f$ et $G_f(\lambda
z)=\log|\lambda|+ G_f(z)$ pour tout $\lambda\in\C^*$. On peut
d\'efinir un courant positif ferm\'e de bidegr\'e $(1,1)$
$T$ de $\P^k$ tel que $\pi^*T=\ddc G_f$. Ce courant s'appelle {\it
courant de Green}. {\it La mesure d'\'equilibre} $\mu$ de $f$ est
d\'efinie par $\mu:=T^k$. C'est la seule mesure de probabilit\'e
invariante par $f$
($f^*\mu=d^k\mu$) qui ne charge par les ensembles pluripolaires.
Plus pr\'ecis\'ement, on a 
$$\lim_{s\rightarrow \infty}\frac{1}{d^{ks}}
\sum_{f^s(z)=a}\delta_z=\mu$$
pour tout $a$ n'appartenant pas \`a un certain ensemble
exceptionnel pluripolaire.
La notation $\delta_z$ d\'esigne la
masse de Dirac en $z$. Le support $J$ de $\mu$ s'appelle {\it
l'ensemble de Julia d'ordre maximal} de $f$. Cet ensemble n'est
pluripolaire dans aucun ouvert qui le rencontre. Si $U$ est un
ouvert rencontrant $J$ alors $\bigcup_{s\geq 1}f^s(U)$ est un
ouvert de compl\'ementaire pluripolaire de $\P^k$.
Briend-Duval ont prouv\'e que 
l'ensemble des points p\'eriodiques r\'epulsifs de
$f$ est dense dans $J$ \cite{BriendDuval1}.
\par
Lorsque $f$ est de plus un endomorphisme polynomial de $\C^k$,
on utilise {\it la fonction de Green} d\'efinie par
$$G(z):=\lim_{s\rightarrow \infty}\frac{1}{d^s}\log^+\|f^s(z)\|$$
o\`u $\log^+(z):=\max(0,\log(z))$.
Cette fonction de Green
est \'egalement continue, plurisousharmonique et \`a
croissance logarithmique quand $z\rightarrow \infty$. Le
courant de Green d\'efini ci-dessus est \'egal \`a $\ddc G$.
L'ensemble de Julia $J$ est contenu dans le bord du compact
$\{G=0\}$. Le compact $\{G=0\}$
s'appelle {\it ensemble de Julia rempli}. Notons
que la fonction de Green, le courant de Green, la mesure
d'\'equilibre et l'ensemble de Julia 
de $f^s$ sont \'egaux \`a ceux de $f$ pour tout $s\geq 1$.
\begin{definition} \rm
On dit qu'au voisinage d'un point $a$
l'ensemble $J$ est {\it lamin\'e par des hypersurfaces
r\'eelles de classe ${\cal C}^\alpha$} s'il
existe un syst\`eme de coordonn\'ees r\'eelles
$(x_1,\ldots, x_{2k})$ de classe ${\cal C}^\alpha$ d'un voisinage
$V$ de $a$ et un ferm\'e $K\subset \R$ tels que $J\cap
V=\{(x_1,\ldots,x_{2k}) \mbox{ avec } x_{2k}\in K\}\cap V$.
On dit que dans un ouvert $U$,
l'ensemble $J$ est
{\it lamin\'e par des hypersurfaces
r\'eelles de classe ${\cal C}^\alpha$} s'il
est au voisinage de tout point de $U$.
\end{definition}
\par
Notre r\'esultat principal est le th\'eor\`eme suivant:
\begin{theoreme} Soit $f$ un endomorphisme polynomial de degr\'e
  $d\geq 2$ de $\C^k$ qui 
  se prolonge en endomorphisme holomorphe
  de $\P^k$ $(k\geq 2)$. Supposons que dans un ouvert de $\C^k$ 
  l'ensemble de Julia d'ordre
  maximal de $f$ soit non vide et soit lamin\'e par des 
  hypersurfaces r\'eelles de classe ${\cal
  C}^{2+\epsilon}$ $(0<\epsilon <1)$. 
Alors pour un syst\`eme de coordonn\'ees convenable de
  $\C^k$, $f$ est homog\`ene. De plus, la restriction de $f$ \`a
  l'hyperplan $\P^k\setminus \C^k \simeq \P^{k-1}$ est une
application de Latt\`es.
\end{theoreme}
\begin{remarque} \rm L'ensemble de Julia d'ordre maximal de $f$
est le bord du bassin d'attraction du point $0$. C'est
une vari\'et\'e r\'eelle analytique, singuli\`ere et
strictement pseudoconvexe \cite{BertelootLoeb1}.
\par 
Au cas d'une variable, les polyn\^omes de degr\'e $d\geq 2$
dont l'ensemble de Julia
est une courbe, sont conjugu\'es \`a $z^d$ ou \`a $\pm
\T_d$ \cite[p. 127]{Steinmetz}
o\`u $\T_d$ est le polyn\^ome de Tchebychev de degr\'e $d$. Le
polyn\^ome $\T_d$ n'est pas homog\`ene.
\end{remarque}
\par
En appliquant le th\'eor\`eme 1.3 \`a un relev\'e polynomial
d'un endomorphisme
 $f$ de $\P^k$ on obtient le r\'esultat suivant:
\begin{corollaire} {\bf \cite{BertelootLoeb2,BertelootLoeb1}} 
Soit $f$ un endomorphisme holomorphe de degr\'e
$d\geq 2$ de $\P^k$
  $(k\geq 1)$. Supposons que dans un ouvert de $\P^k$, le courant de
  Green de $f$ soit lisse et strictement positif. Alors $f$ est une
  application de Latt\`es.
\end{corollaire}
\par
On dit que le courant de Green de $f$ est
{\it lisse et strictement positif}
dans un certain ouvert $U$ si dans $U$ il
est \'egale \`a une forme
de bidegr\'e $(1,1)$ de classe ${\cal C}^{2+\epsilon}$, 
strictement positive et ferm\'ee.
Un petit d\'efaut du lemme 2.5 ne nous permet
pas de remplacer la classe ${\cal C}^{2+\epsilon}$ par ${\cal C}^2$.
\par
Le point de
vu d'Eremenko sur les solutions non triviales de l'\'equation
fonctionnelle $f\circ g=g\circ f$, en particulier sur les fonctions de
Latt\`es, concerne la notion d'orbifold
qui est
consid\'er\'ee par Thurston \cite{Thurston}. Cette notion est
\'egalement utilis\'ee au cas de dimension quelconque.
\begin{definition} {\bf \cite{Dinh2,DinhSibony}} \rm 
Soit $X$ une vari\'et\'e complexe. Notons $\H(X)$
l'ensemble des sous-ensembles analytiques
irr\'eductibles de codimension 1 de $X$.
On appelle {\it orbifold} un
couple $(X,m)$ o\`u $m$ est une
fonction d\'efinie sur $\H(X)$ \`a valeurs dans
$\N^+\cup\{\infty\}$, \'egale \`a 1 sauf sur une
famille localement finie de sous-ensembles analytiques de $X$.
{\it Un
rev\^etement d'orbifolds} $\pi:\ (X_1,m_1)\longrightarrow
(X_2,m_2)$ est un rev\^etement ramifi\'e de $X_1\setminus
\bigcup_{m_1(H)=\infty}H$ dans $X_2\setminus\bigcup_{m_2(H)=
\infty}H$
v\'erifiant $\mult(\pi,H).m_1(H)=m_2(\pi(H))$ pour tout
$H\in\H(X_1)$, o\`u $\mult(\pi,H)$ d\'esigne
la multiplicit\'e de $\pi$
en un point g\'en\'erique de $H$.
\end{definition} 
\par
Soit $f$ un endomorphisme holomorphe de $\P^k$ 
d\'efinissant un rev\^etement d'un orbifold ${\cal
O}=(\P^k,m)$ dans lui-m\^eme. Alors, l'ensemble critique de $f$
est {\it pr\'ep\'eriodique}, c'est-\`a-dire
$f^{n_1}({\cal C}_f)=f^{n_2}({\cal C}_f)$
pour certains $0\leq n_1<n_2$ o\`u ${\cal C}_f$ d\'esigne l'ensemble
critique de $f$. On dit qu'une telle application est {\it
critiquement finie}. La preuve du corollaire suivant se
trouve dans \cite{DinhSibony}:
\begin{corollaire}
Sous l'hypoth\`ese du th\'eor\`eme 1.3 ou du corollaire 1.5, il existe
un orbifold ${\cal O}=(\P^k,m)$ tel que $f$ d\'efinisse un
rev\^etement de ${\cal O}$ dans lui-m\^eme. En particulier, $f$ est
critiquement fini.
\end{corollaire}
\par
Voici le plan de la d\'emonstration.
Nous montrons d'abord qu'il existe un ouvert $U$ dans
lequel $J$ est une
hypersurface ${\cal C}^{2+\epsilon}$ strictement pseudoconvexe. Ceci
s'appuie sur trois principaux arguments:
$G$ ne peut s'annuller au
voisinage d'aucun point de $J$ car $J\subset\partial \{G=0\}$, $J$
n'est pas lamin\'e par des vari\'et\'es
complexes et un th\'eor\`eme de
Tr\'epreau qui implique que si une feuille
$\F$ de la lamination de $J$ ne
contient aucun sous-ensemble analytique
passant par un point $a\in \F$
alors $G$ est nulle dans un demi-voisinage de $a$ bord\'e par $\F$.
\par
On choisit un point p\'eriodique r\'epulsif
$b\in J\cap U$ de p\'eriode
$n$. Posons $g:=f^n$. Le fait que $J$ est
strictement pseudoconvexe et
invariant par $g$ permet de rendre \`a $J$
la forme $\Im z_k=\|z'\|^2$
en utilisant un changement de coordonn\'ees local. Plus pr\'ecis\'ement,
on construit une application holomorphe
(appel\'ee application de Poincar\'e)
$\varphi: \C^k\longrightarrow \C^k$ telle que $J^*:=\varphi(J)=\{\Im
z_k=\|z'\|^2\}$. Les applications biholomorphes locales qui envoient
un ouvert de
$J^*$ dans un autre, en particulier $\Lambda:=\varphi^{-1}\circ g\circ
\varphi$, sont affines. On peut les d\'ecrire. En
suite, on construit le groupe $\A$ de telles applications qui
pr\'eservent les fibres  $\varphi^{-1}(z)$ et on montre que ce groupe
agit transitivement sur tous les fibres de $\varphi$. 
\par
Le dernier paragraphe se consacre \`a la construction de la fibration
invariante qui entra\^{\i}nera l'homog\'enit\'e de $f$. Nous
montrons que l'image de $\{z'=u\}$ par $\varphi$ est contenue dans une
droite complexe pour tout $u\in\C^{k-1}$. Pour
ceci on utilise les propri\'et\'es de la fonction de Green, sa
croissance et sa harmonicit\'e sur les vari\'et\'es stables.
On montre finalement que ces droites se
coupent en un point et cette fibration de droites
est invariante par $f$.
\section{Applications de Poincar\'e}
\ahead
On consid\`ere un endomorphisme
$f$ de degr\'e $d\geq 2$ de
$\C^k$ dont l'ensemble de
Julia d'ordre maximal $J$ est lamin\'e par des hypersurfaces r\'eelles
de classe ${\cal C}^{2+\epsilon}$ dans un certain ouvert $V$
de $\C^k$. 
Nous allons \'etudier $f$ en
utilisant sa forme triangulaire au voisinage des
points p\'eriodiques r\'epulsifs.
\par
Une application $\Lambda$ de $\C^k$ dans lui-m\^eme est appel\'ee 
{\it application triangulaire dilatante} si elle est du type
$$\Lambda(z)=(\lambda_1z_1,\lambda_2 z_2+P_2(z_1),...,\lambda_k
z_k+P_k(z_1,...,z_{k-1}))$$
o\`u $1<|\lambda_1|\leq |\lambda_2|\leq \cdots \leq
|\lambda_k|$ et $P_i$ est un
polyn\^ome en $z_1,...,z_{i-1}$ v\'erifiant
$P_i(\lambda_1z_1,...,\lambda_{i-1}
z_{i-1})= \lambda_iP_i(z_1,...,z_{i-1})$ pour tout $i=2,...,k$.  
\begin{lemme} Soient $f$ un endomorphisme holomorphe de $\C^k$ et
$a$ un point fixe r\'epulsif de $f$. Alors 
  il existe une
  application holomorphe ouverte $\varphi_a: \C^k\longrightarrow \C^k$
  (application de Poincar\'e) et
  une application triangulaire dilatante 
$\Lambda_a :\C^k\longrightarrow \C^k$ telles
  que  $\varphi_a(0)=a$, $\varphi_a'(0)$ inversible et 
  $f\circ \varphi_a =\varphi_a\circ \Lambda_a$. Si de plus $a$
appartient \`a $J$ alors $\varphi(\C^k)$ est un ouvert de
compl\'ementaire pluripolaire de $\P^k$.
\end{lemme}
\begin{preuve} D'apr\`es le th\'eor\`eme de Sternberg
  \cite{Sternberg}, il existe  une
  application holomorphe ouverte $\varphi_a$ d'un voisinage $U$ de $0\in
  \C^k$ dans $\C^k$ et une
  application triangulaire dilatante
$\Lambda_a :\C^k\longrightarrow \C^k$ telles
  que $\varphi_a(0)=a$, $\varphi_a'(0)$ inversible et 
  $f\circ \varphi_a =\varphi_a\circ \Lambda_a$ au voisinage de $0$.
On va prolonger l'application $\varphi_a$ en application
  holomorphe ouverte de $\C^k$ dans $\C^k$.
Comme $\Lambda_a$ est dilatante, 
  pour tout  $z\in\C^k$ fix\'e la suite
  $\Lambda_a^{-m}(z)$ tend vers $0$. 
On pose $\varphi_a(z):=f^m\circ\varphi_a\circ \Lambda_a^{-m}(z)$
  pour $m$ suffisamment grand. Cette application est bien d\'efinie,
  holomorphe et
  ind\'ependante de $m$ car $f^m\circ \varphi_a =\varphi_a\circ
\Lambda_a^m$ au
  voisinage de $0$. De plus, elle est ouverte car $\varphi_a$ est
  ouverte au voisinage de $0$. Par prolongement analytique,
$f\circ \varphi_a=\varphi_a\circ f$ dans tout $\C^k$.
\par
Supposons maintenant que $a\in J$. Par d\'efinition de $\varphi$,
l'ouvert $\varphi(\C^k)$ est \'egal \`a l'ensemble $\bigcup_{m\geq
1}f^m(U)$. Comme $U\cap J\not=\emptyset$, ce dernier ouvert
est de compl\'ementaire pluripolaire ({\it voir} l'introduction). 
\end{preuve}
\begin{corollaire} Soit $f$ un endomorphisme polynomial de $\C^k$ qui
  se prolonge en endomorphisme holomorphe de $\P^k$. Alors son
  ensemble de Julia d'ordre maximal $J$ ne contient aucune vari\'et\'e
complexe passant par un point p\'eriodique r\'epulsif. En
particulier, $J$ n'est pas lamin\'e par des vari\'et\'es
complexes dans aucun ouvert qui le rencontre.
\end{corollaire}
\begin{preuve} Supposons que $J$ contient une vari\'et\'e
complexe $\Sigma$ passant par un point p\'eriodique r\'epulsif $a\in
J$. Quitte \`a remplacer $f$ par un it\'er\'e, on peut supposer
que $a$ est fixe. Notons $\Lambda_a$ et $\varphi_a$ comme dans le
lemme 2.1. Posons $\Sigma^*:=\varphi_a^{-1}(\Sigma)$.
\par
Comme $\Lambda_a$ est dilatante, la famille $\{\Lambda_a^m|_{\Sigma^*}
\mbox{ avec } m\in \N\}$ n'est pas normale. D'apr\`es le lemme de
renormalisation de
Zalcman \cite{DinhSibony,Schwick}, l'adh\'erence de l'ensemble
$\bigcup_{m\in\N} \Lambda_a^m(\Sigma^*)$
contient des images holomorphes non
constantes 
de $\C$. Par cons\'equent, $J^*=\varphi_a^{-1}(J)$ contient des images
holomorphes non constantes de $\C$.
Comme $\varphi_a$ est ouverte, $J$ contient
aussi des images holomorphes non constantes de $\C$. Ceci contredit le
th\'eor\`eme de Liouville car $J$ est un compact de $\C^k$.
\par
D'apr\`es le th\'eor\`eme de Briend-Duval \cite{BriendDuval1},
l'ensemble des points
p\'eriodiques r\'epulsifs est dense dans $J$. On en d\'eduit que
$J$ ne peut \^etre lamin\'e par des vari\'et\'es complexes dans
aucun ouvert qui le rencontre.
\end{preuve}
\begin{proposition} Soit $f$ un endomorphisme de $\C^k$ v\'erifiant
  l'hypoth\`ese du th\'eor\`eme 1.3. Alors 
  il existe un ouvert $U$ de $\C^k$ tel que $J \cap U$
  soit une hypersurface non vide de classe ${\cal C}^{2+\epsilon}$ 
  et strictement
  pseudoconvexe de $U$. 
\end{proposition}
\par
Pour prouver cette proposition, on aura besoin du lemme suivant:
\begin{lemme} Soient $I$ un ouvert de $\R$ et $K$
  un ferm\'e parfait de $I$.
  Alors l'ensemble suivant est dense dans $K$:
  $$Z_K:=\{x\in K,\ ]x-\delta,x[\cap K\not = \emptyset \mbox{ et } 
           ]x,x+\delta[\cap K \not = \emptyset \mbox{ pour tout }
           \delta >0\}.$$
\end{lemme}
\begin{preuve} Supposons que $K$ est non vide. On montre d'abord que 
$Z_K$ est aussi non vide. Supposons que ce n'est pas le cas. Notons
$S:=I\setminus K$. Alors $S$ est une r\'eunion d\'enombrable
d'intervalles disjoints 
$S_n$ de $\R$ et chaque point de $K$ est un sommet de
l'un des $S_n$. En particulier, $K$ est 
d\'enombrable et les composantes connexes de $K$ sont des singletons.
Par hypoth\`ese, $K$ ne contient pas de point
isol\'e. Il est donc un ensemble de Cantor. 
Ceci est impossible car tout ensemble de Cantor est non d\'enombrable.
\par
De m\^eme mani\`ere, on montre que si $I'\subset I$ est un ouvert
v\'erifiant $K':=K\cap I'\not =\emptyset$ alors $Z_{K'}\not
=\emptyset$ et donc $Z_K\cap I'\not=\emptyset$. Ceci implique que $Z_K$
est dense dans $K$. 
\end{preuve}
{\it Preuve de la proposition 2.3 ---}
Montrons d'abord qu'il existe un ouvert $U$ dans lequel $J$ est
une hypersurface de classe ${\cal C}^{2+\epsilon}$.
Par hypoth\`ese sur la lamination de $J$, il existe
un ouvert $V\subset \C^k$ muni des coordonn\'ees locales
$(w_1,\ldots,w_k)$ tel que $J\cap V$ soit une r\'eunion non vide
de graphes de fonctions ${\cal C}^{2+\epsilon}$ au-dessus de
$$B_{2k-1}:=
\{w \mbox{ v\'erifiant } \|w\|<1 \mbox{ et } \Re w_k=0\}.$$  
Notons ces graphes par $\Gamma_\nu$ o\`u le param\`etre $\nu\in \R$
est la partie r\'eelle de la $k$-i\`eme coordonn\'ee du point
$\Gamma_\nu\cap\{w'=0 \mbox{ et } \Im w_k=0\}$.
C'est-\`a-dire $\Gamma_\nu\cap
\{w'=0 \mbox{ et } \Im z_k=0\}=(0,\nu)$.
Notons $K$ l'ensemble des valeurs de $\nu$.
C'est un ferm\'e d'un certain ouvert $I$ de $\R$. 
\par
Supposons que $J$ n'est une hypersurface dans aucun ouvert qui le
rencontre. Alors $K$ est un ensemble parfait.  L'ensemble
$Z_K$ est d\'efini dans le lemme 2.4.
Fixons un $\nu\in Z_K$.
Montrons que $\Gamma_\nu$ est Levi-plate. Supposons que
$\Gamma_\nu$ n'est pas Levi-plate. D'apr\`es le th\'eor\`eme de
Tr\'epreau, il existe un point
$w\in\Gamma_\nu$ au voisinage duquel les disques holomorphes
attach\'es \`a
$\Gamma_\nu$ remplit un demi-voisinage \`a un c\^ot\'e bord\'e par
$\Gamma_\nu$ \cite{Trepreau}. Notons $W$ ce demi-voisinage. La
fonction de Green $G$ est plurisousharmonique,
positive ou nulle. Elle est \'egale \`a $0$ sur
$J$. Par principe du maximum, 
elle s'annulle sur $W$. Par d\'efinition de $K$ et de $Z_K$, il
existe des feuilles de $J$ qui rencontrent $W$. On en d\'eduit
qu'il existe des points de $J$ au voisinage desquels $G$ est nulle.
Ceci contredit le fait que $J$ est contenu dans le bord de $\{G=0\}$.
\par
Alors $\Gamma_\nu$ est Levi-plate pour tout $\nu\in Z_K$.
D'apr\`es le lemme 2.4, $Z_K$ est dense dans $K$. Par
continuit\'e, 
$\Gamma_\nu$ est Levi-plate pour tout $\nu\in K$.
En particulier, $J$ est lamin\'e par des vari\'et\'es
complexes. Ceci contredit le corollaire 2.2.
\par
On conclut qu'il existe un ouvert $U$ de $\C^k$
dans lequel $J$ est une
vari\'et\'e de classe ${\cal C}^{2+\epsilon}$. Comme $G$ est
plurisouharmonique et $J$ est contenu dans le bord de $\{G=0\}$,
la vari\'et\'e $J\cap U$ est
pseudoconvexe en tout point. Montrons qu'il est strictement
pseudoconvexe en au moins un point. Supposons que ceci n'est pas
le cas. Pour tout $a\in J\cap U$, on note $L_a$ l'espace tangent
complexe de $J$ en $a$, $\L_a$ la forme de Levi de $J$ en $a$
et $K_a\subset L_a$ l'ensemble des
vecteurs holomorphes $X$ v\'erifiant $\L_a(X,\overline X)=0$.
Comme $J$ est pseudoconvexe, $\L_a$ est d\'efinie par une matrice
semi-positive diagonable.
Par cons\'equent,  $\L_a(X,\overline X)=0$
entra\^{\i}ne  
$\L_a(X,.)=0$. Alors $K_a=\{X\in L_a, \L_a(X,.)=0\}$ 
est un espace vectoriel complexe de dimension
au moins 1. Quitte \`a remplacer $U$ par un ouvert convenable, on
peut supposer que la dimension $r$ de $K_a$ est ind\'ependante de
$a \in J\cap U$.
D'apr\`es le th\'eor\`eme de Freeman \cite[p.185]{Boggess}, $J\cap
U$ est lamin\'e par des vari\'et\'es complexes de dimension $r$.
Ceci contredit le corollaire 2.2.
\par
Alors $J\cap U$ est strictement pseudoconvexe en au moins un
point. Quitte \`a remplacer $U$ par un petit voisinage de ce
point, on peut supposer que $J\cap U$ est une hypersurface
strictement pseudoconvexe de $U$.
\par
\hfill $\square$
\par
Par suite, on consid\`ere 
$b\in J\cap U$ un point p\'eriodique r\'epulsif de p\'eriode $n$ de
$f$. Posons $g:=f^n$. D'apr\`es le lemme 2.1, il existe une
application  holomorphe ouverte $\varphi: \C^k\longrightarrow \C^k$
et une
application triangulaire dilatante
$\Lambda :\C^k\longrightarrow \C^k$ telles
que $\varphi(0)=b$, $\varphi'(0)$ inversible et
$g\circ \varphi =\varphi\circ \Lambda$.
Posons $J^*:=\varphi^{-1}(J)$, $\mu^*:=\varphi^*(\mu)$ et
$G^*:=G\circ\varphi$.  L'application triangulaire dilatante
$\Lambda=(\Lambda_1,\ldots,\Lambda_k)$
s'\'ecrit sous la forme
$$\Lambda(z)=(\lambda_1 z_1, \lambda_2 z_2
+P_2(z_1),...,\lambda_k z_k + P_k(z_1,\ldots,z_{k-1}))$$
o\`u $1<|\lambda_1|\leq |\lambda_2|\leq \cdots \leq |\lambda_k|$ et
$P_i$ est un polyn\^ome en $z_1,...,z_{i-1}$ v\'erifiant
$$P_i(\lambda_1z_1,...,\lambda_{i-1}z_{i-1})=\lambda_i
P_i(z_1,\ldots, z_{i-1}).$$
\begin{lemme}  On a 
$\lambda_k=|\lambda_i|^2>1$ pour tout $i=1,...,k-1$. 
La vari\'et\'e $J^*$ est
d\'efinie par $\Im (\alpha z_k)
= P(z',\overline z')$ o\`u $\alpha\not =0$
est un certain nombre complexe et
$P$ est un polyn\^ome r\'eel
homog\`ene de degr\'e $2$.
\end{lemme}
\begin{preuve} 
On choisit un syst\`eme lin\'eaire de coordonn\'ees
$(w_1,\ldots,w_k)$ de $\C^k$ de sorte que l'espace tangent de
$J^*$ en $0$ soit \'egal \`a $H:=\{w'=0 \mbox{ et } \Im w_k=0\}$. 
Dans ces coordonn\'ees $\Lambda$ n'est pas forc\'ement
triangulaire. On note $\Lambda^*_1,...,\Lambda^*_k$ leurs fonctions
coordonn\'ees.
Au voisinage de $0$, $J^*$ est le graphe d'une fonction r\'eelle
$h$ de classe
${\cal C}^{2+\epsilon}$ qui est d\'efinie sur un
ouvert de $H$, c'est-\`a-dire $J^*=\{\Im w_k=h(w',\Re
w_k)\}$ au voisinage de $0$. Consid\'erons le d\'eveloppement
d'ordre 2 de $h$:
$$h(w',\Re w_k)=P(w')+S(w') \Re w_k+ c(\Re w_k)^2+
\O(\|(w',\Re w_k)\|^{2+\epsilon})$$
o\`u $P$ (resp. $S$) est un polyn\^ome r\'eel homog\`ene
de degr\'e 2 (resp. 1) et $c$ est une constante.
Comme $J^*$ est strictement pseudoconvexe en $0$, la partie
non-harmonique $Q$ dans $P$ est d\'efinie strictement positive.
\par
Le fait que $J^*$ est invariant par $\Lambda$ implique que
l'espace tangent complexe $\{w_k=0\}$ et l'espace tangent r\'eel
$\{\Im w_k=0\}$ de $J^*$ en $0$ sont
invariants par la d\'eriv\'ee
$\Lambda'(0)$. C'est-\`a-dire il existe $\theta\in \R$ tel que
$\Lambda^*_k(w)=\theta w_k+\O(\|w\|^2)$. On \'ecrit
$(\Lambda^*_1,...,\Lambda^*_{k-1})$ sous la forme
$$(\Lambda^*_1,...,\Lambda^*_{k-1})=Aw'+Bw_k+\O(\|w\|^2)$$
o\`u $A$ est une matrice carr\'ee
de taille $k-1$ et $B\in \C^{k-1}$. On obtient du d\'eveloppement
de $h$ et de l'invariance de $J^*$ par $\Lambda$ en consid\'erant
seulement les parties non-harmoniques d'ordre 2 en $w'$ que 
$\theta Q(w') = Q(Aw')$.
Soient $\nu$ une valeur propre de $A$ et
 $u\not =0$ un vecteur propre associ\'e \`a $\nu$. On a $\theta
 Q(\lambda u)=Q(\lambda \nu u)$ pour tout $\lambda\in \C$. 
Comme $Q$ est d\'efinie strictement positive,
on obtient de la derni\`ere \'egalit\'e que
 $\theta=|\nu|^2$. 
\par
L'application $\Lambda$ \'etant dilatante, le module de $\nu$ doit
\^etre strictement plus grand \`a 1. 
On constate 
que $\theta$ est la plus grande valeur
propre de $\Lambda'(0)$ et elle est strictement plus grande que les
autres en module. Ceci implique que $\theta=\lambda_k$ et
$\{z_k=0\}=\{w_k=0\}$. On peut
alors choisir $w'=z'$ et $w_k=\alpha z_k$ avec un $\alpha\in \C^*$
convenable. On a $\lambda_k=|\lambda_i|^2$
pour tout $i=1,...,k-1$. L'application $\Lambda$ \'etant
triangulaire, les polyn\^omes $P_i$ sont
de degr\'e 1 si $i\not=k$ et de degr\'e 2 si $i=k$. Comme
$\Lambda$ est triangulaire, on a $B=0$.
Pour simplifier les notations, on suppose que
$\alpha=1$. Le fait que $J^*$ est invariant par $\Lambda$
implique que
$$\lambda_k h(z',\Re z_k) +\Im P_k(z') =
h(\Lambda_1(z'),...,\Lambda_{k-1}(z'),\lambda_k
\Re z_k +\Re P_k(z')).$$
En consid\'erant les termes d'ordre 2, on obtient:
$$\lambda_k[P(z')+S(z')\Re z_k +c(\Re z_k)^2]+\Im P_k(z')
=P(Az')+\lambda_k S(Az')\Re z_k +
c\lambda_k^2(\Re z_k)^2.$$ 
Cette derni\`ere \'equation
implique que $S=0$ et $c=0$
car les valeurs propres de $A$ est de module
strictement sup\'erieur \`a 1 et $\lambda_k>1$.
\par
Posons $r(z',\Re z_k):=h(z',\Re z_k)-P(z')
=\O(\|(z',\Re z_k)\|^{2+\epsilon})$.
Il reste \`a montrer que $r=0$.
Soit $\Psi:\C^{k-1}\times
\R\longrightarrow \C^{k-1}\times \R$ d\'efinie
par
$$\Psi(z',\Re z_k):=(\Lambda_1(z'),...,\Lambda_{k-1}(z'),
\lambda_k \Re z_k+
\Re P_k(z')).$$
Alors on a $\lambda_k r(z',\Re z_k)=r(\Psi(z',\Re z_k))$. On en
d\'eduit que $r(z',\Re z_k)=\lambda_k^s r(\Psi^{-s}(z',\Re 
z_k))$. D'autre part, on obtient de la triangularit\'e de
$\Lambda$ que  $\Psi^{-s}(z',\Re
z_k)=\O({\lambda_k}^{-s/2}\log s)$ quand
$s\rightarrow \infty$. Par suite,
$$r(z',\Re z_k)=\lim_{s\rightarrow \infty} \lambda_k^s
r(\Psi^{-s}(z',\Re z_k))=\lim_{s\rightarrow \infty}
\lambda_k^s \O(\lambda_k^{-(2+\epsilon)s/2}(\log
s)^{2+\epsilon})=0.$$ 
\end{preuve}
\begin{proposition} Quitte \`a effectuer un changement de
coordonn\'ees qui fixe la droite $\{z'=0\}$, on
  peut supposer que $J^*$ soit d\'efinie par l'\'equation 
  $\Im z_k=\|z'\|^2$. Dans ces nouvelles coordonn\'ees
$\Lambda$ n'est pas forc\'ement triangulaire et la
fonction $G^*(z)$ est nulle si $\Im z_k\geq \|z'\|^2$. 
\end{proposition}
\begin{preuve} On effectue le premier changement de coordonn\'ees
  $z_k\mapsto \alpha z_k$ qui permet de supposer que $\alpha=1$.
\par
Le polyn\^ome r\'eel $P(z',\overline z')$, qui d\'efinit la
  vari\'et\'e $J^*$, est 
homog\`ene de degr\'e 2. On peut l'\'ecrire en somme
$q(z',\overline z') + \Im H(z')$ o\`u $q$ est une forme
hermitienne et
$H$ est un polyn\^ome holomorphe de degr\'e 2. Quitte \`a changer
$z_k$ par $z_k-H(z')$, on peut supposer que $H=0$. 
\par
Comme $J$ est strictement pseudoconvexe en $b$, $q$ est d\'efinie
strictement positive ou n\'egative. Quitte \`a remplacer $z_k$ par
$-z_k$ au cas n\'ecessaire, on peut supposer que $q$ est d\'efinie
strictement positive. Un changement lin\'eaire
convenable des coordonn\'ees
$z_1,...,z_{k-1}$ rend $q$ sa forme canonique $q(z',\overline
z')=\|z'\|^2$.
\par
Il est clair que les changements de coordonn\'ees
effectu\'es fixent la droite
$\{z'=0\}$. La fonction
plurisouharmonique, positive $G^*$ est nulle sur l'hypersurface
$J^*$ qui est strictement pseudoconvexe.
Par principe du maximum, elle doit s'annuller \'egalement sur
les disques holomorphes \`a bord dans $J^*$. On en d\'eduit
qu'elle s'annulle sur 
$\{z\in\C^k, \Im z_k\geq \|z'\|^2\}$. 
\end{preuve}
\section{Groupes d'automorphismes pr\'eservant $J^*$}
\ahead
On utilise les notations du paragraphe 2 avec
le syst\`eme de
coordonn\'ees d\'efini par la proposition 2.6.
On n'aurra plus besoin du fait que $\Lambda$ est
triangulaire.
Rappelons que la droite $\{z'=0\}$ est toujours invariante par
$\Lambda$ et la restriction de $\Lambda$ sur cette droite est 
l'application lin\'eaire $(0,z_k)\mapsto (0,\lambda_k z_k)$ o\`u 
$\lambda_k\in\R^+$ est la plus grande valeur propre de
$\Lambda'(0)$. Les autres valeurs propres de $\Lambda'(0)$ sont
de module $\sqrt{\lambda_k}$.
\par
Pour tout $w=(w',v+i\|w\|^2)\in J^*$ avec $w'\in\C^{k-1}$ et
$v\in\R$, posons $\tau_w(z):=(z'+w',
z_k+2i\overline w' z'+v+i\|w'\|^2)$.
C'est une application affine,
inversible et v\'erifiant $\tau_w(0)=w$. Elle pr\'eserve les
vari\'et\'es $\{\Im z_k-\|z'\|^2=\const\}$.
Notons
pour tout $\lambda$ r\'eel positif $\sigma_\lambda$
l'application lin\'eaire
de $\C^k$ dans lui-m\^eme d\'efinie par
$\sigma_\lambda(z',z_k):=(\sqrt{\lambda} z',
\lambda z_k)$. Cette application 
pr\'eserve $J^*$. Notons $\G_1$ (resp. $\G_2$, $\G_3$ et $\G_4$) 
l'ensemble des automorphismes affines
de $\C^k$ qui s'\'ecrivent sous forme $(u(z'),z_k)$
(resp. $\sigma_\lambda\circ (u(z'),z_k)$,
$\tau_w\circ (u(z'),z_k)$ et $\tau_w\circ \sigma_\lambda \circ
(u(z'),z_k)$) o\`u 
$u$ est un automorphisme lin\'eaire isom\'etrique de $\C^{k-1}$,
$\lambda\in\R^+$ et $w\in J^*$. On v\'erifie sans peine que $\G_1$
et $\G_2$ sont des groupes. L'ensemble
$\G_3$ est \'egal \`a l'ensemble des
\'el\'ements de $\G_4$ dont toute valeur propre de la partie
lin\'eaire est de module 1.
\begin{lemme} Soient $w$, $\tilde w$
deux points de $J^*$ et $\tau$ une
  application holomorphe ouverte d'un voisinage
$W$ de $w$ dans $\C^k$ avec
  $\tau(w)=\tilde w$.
Supposons que $\tau(J^*\cap W)\subset J^*$. Alors 
 $\tau_{\tilde w}^{-1}\circ\tau\circ\tau_w$ appartient \`a
$\G_2$. Par cons\'equent, $\G_3$ et $\G_4$ sont des groupes et
$\tau\in \G_4$.
\end{lemme}
\begin{preuve} On peut consid\'erer ce lemme comme un corollaire
d'un r\'esultat d'Alexander \cite{Alexander}. Donnons ici une
autre preuve.
\par
Observons que
l'application $\tau_{\tilde w}^{-1}\circ \tau\circ \tau_w$
  pr\'eserve $J^*$ et fixe le point 0. Il suffit donc de consid\'erer
  le cas $w=\tilde w=0$. Dans ce cas $\tau_w=\tau_{\tilde w}=\Id$.
On pose
  $\tau=(\tau',\tau'')$ o\`u $\tau'$ se compose par $k-1$
premi\`eres fonctions coordonn\'ees de $\tau$ et $\tau''$ est sa
derni\`ere fonction coordonn\'ee.
\par
Le fait que $J^*$ est invariant par $\tau$ implique que
$$\Im\tau''(z',\Re z_k+i\|z'\|^2)=\|\tau'(z',\Re
z_k+i\|z'\|^2)\|^2.$$
En particulier, pour $\Re z_k=0$ on a
$$\Im\tau''(z',i\|z'\|^2)=\|\tau'(z',i\|z'\|^2)\|^2.$$
La partie pluriharmonique du membre \`a droite est nulle. Par
cons\'equent, $\Im\tau''(z',0)=0$ et donc $\tau''(z',0)=0$. On
peut dire que $\tau$ pr\'eserve l'hyperplan $\{z_k=0\}$.
\par
Soient $z\in J^*$ assez proche de $0$ et $\tilde z:=\tau(z)$. Posons
$\tilde
\tau=(\tilde\tau',\tilde\tau''):=\tau_{\tilde z}^{-1}
\circ\tau\circ\tau_z$. 
Comme $\tau$, l'application $\tilde \tau$
pr\'eserve $J^*$ et fixe $0$. On
montre de m\^eme mani\`ere que $\tilde \tau$ pr\'eserve
$\{z_k=0\}$. Ceci implique que l'image de l'hyperplan tangent de
$J^*$ en $z$ par $\tau$ est un hyperplan. Soit $\P^{k*}\simeq \P^k$
l'ensemble des
hyperplans de $\P^k$. 
Notons $\F$ la famille des hyperplans dont l'image par $\tau$ est
un hyperplan. C'est un ensemble analytique dans $\P^{k*}$
qui contient la famille $\F'$ des
hyperplans tangents \`a $J^*$ au voisinage de $0$.
Comme $\F'$ est de codimension r\'eelle 1 dans $\P^{k*}$,
la famille analytique $\F$ est un ouvert de $\P^{k*}$ et
contient tous les hyperplans qui passent au voisinage de $0$. On en
d\'eduit que $\tau$ est affine et donc lin\'eaire car $\tau(0)=0$.
\par
Comme $J^*$ est invariante par $\tau$ et $\tau(0)=0$, l'hyperplan
tangent r\'eel $\{\Im z_k=0\}$ de $J^*$ en $0$ est invariant par
$\tau'(0)$.
Par cons\'equent, 
$\frac{\partial \tau''}{\partial z_k}(0)$ est un nombre r\'eel et
$\frac{\partial \tau''}{\partial z'}(0)=0$. On peut donc \'ecrire
$\tau''(z)=\lambda z_k$
o\`u $\lambda\in\R$.
Le fait que $J^*$ est invariant par $\tau$ implique
que $\lambda \|z'\|^2= \|\tau'(z',\Re z_k+i\|z'\|^2)\|^2$.
D'o\`u $\lambda>0$, $\tau'$ est ind\'ependant de $z_k$ et
$\|\frac{1}{\sqrt{\lambda}}\tau'(z')\|=\|z'\|^2$. Ceci implique que
$u:=\frac{1}{\sqrt{\lambda}}\tau'(z')$ est une isom\'etrie complexe
de $\C^{k-1}$. L'application $\tau(z)$ s'\'ecrit sous la forme
$\sigma_\lambda\circ
(u(z'),z_k)$. Il appartient donc \`a $\G_2$.
\par
Notons $\tilde \G_4$ le groupe engendr\'e par les \'el\'ements de
$\G_4$. Soit $\tau^*\in\tilde \G_4$. On montre que $\tau^*\in \G_4$.
Comme $\tau^*$ est un automorphisme pr\'eservant $J^*$, le point
$w^*:=\tau^*(0)$ appartient \`a $J^*$. D'apr\`es la premi\`ere
partie, $\tau_{w^*}^{-1}\circ\tau^*$ appartient \`a $\G_2$. Ceci
implique que $\tau^*\in \G_4$. L'ensemble $\G_4$ est donc un groupe.
En particulier, $\tau\in \G_4$.
\par
L'ensemble $\G_3$ est \'egal \`a l'ensemble des \'el\'ements de
$\G_4$ dont toute valeur propre de la partie lin\'eaire
est de module 1. Par
cons\'equent, $\G_3$ est aussi un groupe.
\end{preuve}
\begin{corollaire} L'application $\Lambda$ s'\'ecrit sous la
forme $\sigma_{\lambda_k}\circ \tau $ o\`u $\tau$ est un \'el\'ement
de $\G_1$.
\end{corollaire}
\begin{preuve} D'apr\`es le lemme 3.1, l'application $\Lambda$
s'\'ecrit sous la forme $\sigma_\lambda\circ \tau$ o\`u
$\lambda>0$ et $\tau\in \G_1$. On
sait que la restriction de $\Lambda$ sur $\{z'=0\}$ est \'egale
\`a $(0,z_k)\mapsto (0,\lambda_k z_k)$. Par cons\'equent,
$\lambda=\lambda_k$ et donc $\Lambda(z)=\sigma_{\lambda_k}\circ
\tau$.
\end{preuve}
\begin{corollaire} Soient $w$, $\tilde w$
deux points de $J^*$ v\'erifiant
  $\varphi(w)=\varphi(\tilde w)$. Alors il existe 
  $\tau\in \G_3$ avec $\tau(w)=\tilde w$ tel que
  $\tau\circ\varphi =\varphi$. 
\end{corollaire}
\begin{preuve} On consid\`ere le cas o\`u $w$ et $\tilde w$
n'appartiennent pas \`a l'ensemble critique de $\varphi$. Le cas
g\'en\'eral sera d\'eduit par prolongement analytique. Fixons des
petits voisinages $W$, $\tilde W$, $V$ de $w$, $\tilde w$ et de
$\varphi(w)=\varphi(\tilde w)$ tels que $\varphi$ r\'ealise un
biholomorphisme entre $W$ (resp. $\tilde W$) et $V$. On pose
$\tau:=(\varphi|_{\tilde W}^{-1})\circ \varphi|_W$. On a
$\varphi\circ\tau=\varphi$ sur $W$ et 
$$\tau(J^*\cap W)=(\varphi|_{\tilde W})^{-1}(J\cap V)=J^*\cap
\tilde W.$$   
D'apr\`es le lemme 3.1, l'application $\tau$ appartient \`a $\G_4$. 
Par prolongement analytique on a
$\varphi\circ\tau=\varphi$ sur $\C^k$.
\par
L'application $\tau$ s'\'ecrit sous forme $\tau^* \circ 
\sigma_\lambda\circ \tilde \tau$ o\`u $\lambda >0$, $z\in J^*$ et
$\tau^*$ et $\tilde \tau$ appartient \`a $\G_3$.
Il reste \`a montrer que $\lambda=1$. Supposons
que $\lambda\not=1$.
Alors les valeurs propres de la partie lin\'eaire de $\tau$ sont
diff\'erentes de 1. 
Il existe donc $z\in\C^k$ tel que $\tau(z)=z$. Au
voisinage de $z$ l'application $\tau$ r\'ealise donc une permutation
dans les fibres de $\varphi$. Par cons\'equent, pour un certain $m\geq
1$, on a $\tau^m=\Id$. En particulier, les valeurs propres de la
partie lin\'eaire de $\tau$ sont de module 1.
Ceci contredit
le fait que $\lambda\not =1$.
\end{preuve}
\par
Il est clair que
l'ensemble des applications v\'erifiant le corollaire 3.3 forme un
sous-groupe discret de $\G_3$. 
Notons $\A$ ce groupe. Les \'el\'ements de ce groupe
pr\'eservent $J^*$. On note $\A|_{J^*}$ la restriction de $\A$ 
sur $J^*$ ({\it voir} les premi\`eres notations de l'article).
\begin{proposition} Le groupe $\A|_{J^*}$  est co-compact.
Le groupe $\A$ agit transitivement sur les fibres de $\varphi$.
\end{proposition}
\begin{preuve}
Pour ne pas confondre les notations, on note $w'$ les $k-1$
premi\`eres coordonn\'ees de $w$ et $w''$ sa derni\`ere
coordonn\'ee pour tout $w\in\C^k$. Pour toute application $\tau$
de $\C^k$ dans $\C^k$, $\tau'$ se compose par les $k-1$ premi\`eres
fonctions coordonn\'ees de $\tau$ et $\tau''$ est la derni\`ere
fonction coordonn\'ee de $\tau$.
\par
L'application $\Lambda$ \'etant dilatante et pr\'eservant les
fibres de $\varphi$, 
pour montrer que $\A_{|J^*}$ est co-compact, il
  suffit de montrer qu'il existe un ouvert born\'e $W$ de $J^*$
  contenant $0$ tel que pour tout $w\in \partial W$,
il existe $w^*\in
  W$ v\'erifiant $\varphi(w^*)=\varphi(w)$. En effet, ceci implique
  que $W$ contient un domaine fondamental de $\A|_{J^*}$.
\par  
Notons $\Pi$ la projection de $J^*$ dans $\C^{k-1}\times
  \R$ d\'efinie par $\Pi(z',z_k):=(z',\Re z_k)$. On
choisit un nombre fini de points
  $a_s=(a_s',a_s'')\in J^*$ pour $s=1$, 2, 3, ... v\'erifiant les
  propri\'et\'es suivantes:
\begin{enumerate}
\item Pour tout $s$, on a 
$1/2<\|\Pi(a_s)\|<1$.
\item Dans $\C^{k-1}\times \R$,
l'enveloppe convexe $U$ des points $\Pi(a_s)$ contient la boule de
centre $0$ est de rayon $1/2$. 
\item Pour tout $u\in \partial U$, il existe un $s$ tel que
  $\|u-\Pi(a_s)\|<1/16$.  
\end{enumerate}
On choisit $b_0\in J^*$ tel que $\|b_0\|<1/100$ et tel que la mesure
$$\frac{1}{(d_g)^{sk}}\sum_{g^s(z)=\varphi(b_0)}\delta_z$$
tende vers $\mu$ quand $s$ tend vers l'infini. Alors l'ensemble
$\bigcup_{s\geq 0} g^{-s}(\varphi(b_0))$ est dense dans $J$. 
Quitte \`a pertuber l\'eg\`erement les points $a_s$, on peut supposer
qu'il existe un $p$ tel que $g^p(\varphi(a_s))=\varphi(b_0)$ pour tout
$s$. Posons $b_s:=\Lambda^p(a_s)$. Soit $V$ l'enveloppe convexe des
$\Pi(b_s)$ et $W:=\Pi^{-1}(V)$. Observons que $V$ est de taille
$\lambda_k^{p/2}$
dans les directions complexes
et $\lambda_k^p$ dans la direction r\'eelle ({\it voir} le
corollaire 3.2).
\par
Pour $w\in \partial W$, on choisit un $s$
tel que $\|\Pi(\Lambda^{-p}(w))-\Pi(a_s)\|<1/16$.
Si $\tau=(\tau',\tau'')\in
\A$ est tel que $\tau(b_0)=b_s$ et 
$w^*:=\tau^{-1}(w)$ alors $\varphi(w^*)=\varphi(w)$. On montre
que $w^*\in W$. D'apr\` es la condition 2, il suffit de montrer que
$\|{w^*}'\|<\lambda_k^{p/2}/4$ et $\|\Re{w^*}''\|<\lambda_k^p/4$. 
Comme $\tau'$ et $\frac{1}{\sqrt{\lambda_k}}\Lambda'$
sont des isom\'etries de
$\C^{k-1}$, 
on obtient 
$$\|{w^*}'-b_0'\|=\|w'-b_s'\|=\lambda_k^{p/2}
\|\Lambda^{-p}(w)'-a_s'\|<\lambda_k^{p/2}/16.$$
D'o\`u $\|{w^*}'\|<\lambda_k^{p/2}/4$ car $\|b_0\|<1/100$.
\par
Posons $\tilde \tau=(\tilde
\tau',\tilde \tau''):=
\tau_{b_s}^{-1}\circ\tau\circ\tau_{b_0}$,
$w_1:=\tau_{b_s}^{-1}(w)$ et $w_2:=\tilde\tau^{-1}(w_1)$. 
Alors $w^*=\tau_{b_0}(w_2)$ et $\tilde \tau$ est
lin\'eaire isom\'etrique car $\tau\in\G_3$ et $\tau(0)=0$. On a 
$$\|w_2'\|=\|w_1'\|=\|w'-b_s'\|\leq
\lambda_k^{p/2}/16.$$
On a aussi $\Re w_1''=w''-\Re b_s''+2\Im (w'-b_s') (\overline
b_s')$. D'o\`u 
$$\|\Re w_1''\|\leq \|w''-\Re b_s''\|+
2\|w'-b_s'\|\|b_s'\|\leq 3\lambda_k^p/16$$
car
$$\|w''-\Re b_s''\|=\lambda_k^p\|\Lambda^{-p}(w)''
-a_s''\|<\lambda_k^p/16$$
et
$$\|w'-b_s'\|\|b_s'\|=\lambda_k^p \|\Lambda^{-p}(w)'-a_s'\|\|a_s'\|
\leq \lambda_k^p/16.$$
Comme $\tilde \tau$ est lin\'eaire isom\'etrique, on a
$\|\Re w_2''\|\leq 3\lambda_k^p/16$ car $w_2=\tilde \tau^{-1}(w_1)$.
Comme $w^*=\tau_{b_0}(w_2)$, on a $\Re {w^*}''=\Re
w_2''+\Re b_0-2\Im (w_2' \overline b_0)$. Par cons\'equent,
$$\|\Re {w^*}''\|\leq \|\Re w_2''\|+\|\Re b_0''\| +
2\|w_2'\|\|b_0\|\leq
\frac{3\lambda_k^p}{16}+\frac{1}{100}+\frac{\lambda_k^{p/2}}{800}
\leq \frac{\lambda_k^p}{4}.$$
\par
Il en r\'esulte que $\A_{|J^*}$ est co-compact. L'ensemble
$\varphi(J^*)$ est un compact de $J$. D'apr\`es le lemme 2.1,
$J\setminus \varphi(J^*)$ est un ouvert pluripolaire de $J$. Cet
ouvert est donc vide car $J$ n'est pas pluripolaire dans aucun
ouvert qui le rencontre. On a $J=\varphi(J^*)$.
\par
Montrons que $\A$ agit transitivement sur les fibres de $\varphi$.
Le corollaire 3.3 montre que $\A$ agit
transitivement sur le fibre $\varphi^{-1}(z)$ pour tout $z\in J$.
Comme $\# g^{-1}(z)\leq (d_g)^k$ pour tout $z$,
l'ensemble $\Lambda^{-1}(\A w)$ est une
r\'eunion d'au plus $(d_g)^k$ orbites de $\A$
pour tout  $w\in J^*$. Si $\varphi(w)\not\in g({\cal C}_g)$, on a
$\# g^{-1}(\varphi(w))=(d_g)^k$ o\`u ${\cal C}$ signifie
l'ensemble critique.
La relation
$g\circ\varphi=\varphi\circ\Lambda$ implique que $\varphi^{-1}({\cal
C}_g)\subset \Lambda^{-1}({\cal C}_\varphi)$. 
Alors pour tout $w\in J^*\setminus
\Lambda^{-1}({\cal C}_\varphi)$,
l'ensemble $\Lambda^{-1}(\A w)$ est une
r\'eunion de $(d_g)^k$ orbites de $\A$. 
Pour un $w\in J^*\setminus \Lambda^{-1}({\cal C}_\varphi)$ 
fix\'e, on peut choisir
  $\tau_1$, ..., $\tau_{(d_g)^k}$ des \'el\'ements de $\A$
  avec $\tau_1=\Id$ 
  tels que $\Lambda^{-1}(\A w)$ soit la r\'eunion des
orbites $\A.\Lambda^{-1}(\tau_s(w))$ de $\A$ pour $s=1,...,(d_g)^k$.
\par
Par prolongement analytique, ceci est vrai pour tout $w\in \C^k$ \`a
l'exception d'une hypersurface complexe qu'on notera $\Sigma$ o\`u les
points $\varphi(\Lambda^{-1}(\tau_s(w)))$ ne sont pas distincts.
Supposons
  que $\A$ n'agisse pas transitivement sur les fibres de
  $\varphi$. Il existe donc $w$ et $u$ tels que $w\not\in \A.u$
  et $\varphi(w)=\varphi(u)$. 
Quitte \`a pertuber l\'eg\`erement $w$ et $u$ on peut supposer que
$\Lambda^{-n}(\A w)\cap \Sigma =\emptyset$
et $\Lambda^{-n}(\A u)\cap \Sigma
=\emptyset$ pour tout $n\geq 1$.
Par suite,
$\varphi(\Lambda^{-1}(\tau_s(w))$ (resp.
$\varphi(\Lambda^{-1}(\tau_s(u))$) sont les $(d_g)^k$ pr\'eimages
diff\'erentes de $\varphi(w)=\varphi(u)$ par $g$. On en d\'eduit
qu'il existe $s_1$ tel que
$\varphi(\Lambda^{-1}\circ\tau_1(w))=
\varphi(\Lambda^{-1}\circ\tau_{s_1}(w))$. Comme
$\tau_1=\Id$, on a $\varphi(\Lambda^{-1}(w))=
\varphi(\Lambda^{-1}(\tau_{s_1}(w)))$.
\par
Par r\'ecurrence, il existe $s_1$, $s_2$, ... tels que pour tout
$m\geq 1$ on ait
$\varphi(w_m)=\varphi(u_m)$ o\`u $w_m:=\Lambda^{-m}(w)$ et
$u_m:=\Lambda^{-1}\circ
\tau_{s_m}\circ\cdots\circ\Lambda^{-1}\circ\tau_{s_1}(u)$. Comme
$w\not\in \A.u$, on a $w_m\not \in\A.u_m$ pour tout $m\geq 1$.
La suite $u_m$ est born\'ee. En effet $\Lambda^{-1}$ est contractante
et le corollaire 3.3 permet de
contr\^oler les $\tau_{s_j}$. Soit $(u_{m_r})$ une sous-suite
de $(u_m)$ tendant vers
$u_0$. Comme $w_m$ tend vers
$0$, on a  $\varphi(u_0)=\varphi(0)$. Comme $\varphi(0)=b\in J$,
on a $u_0\in J^*$. Il existe donc $\tau\in \A$
tel que $\tau(u_0)=0$. La suite
$v_{m_r}:=\tau(u_{m_r})$
tend vers $0$ et v\'erifie
$\varphi(v_{m_r})=\varphi(w_{m_r})$, $v_{m_r}\not=w_{m_r}$ car
$w_{m_r}\not\in \A.u_{m_r}$.
C'est la contradiction recherch\'ee car $\varphi$ est injective
au voisinage de $0$.
\end{preuve}
\begin{corollaire} L'ensemble $J$ est \'egal \`a $\varphi(J^*)$.
De plus,
pour tout $z\in J$ l'ensemble $\bigcup_{n\geq 0}g^{-n}(z)$ est
dense dans $J$. On en d\'eduit que pour tout $w\in J^*$
l'ensemble
$$\B_w:= \{\Lambda^{-n}\circ\tau (w) \mbox{ avec } n\geq 0 \mbox{ et }
\tau\in \A\}$$
est dense dans $J^*$.
\end{corollaire}
\begin{preuve} On a montr\'e dans la preuve de la proposition 3.4
que $J=\varphi(J^*)$.
On en d\'eduit que $J$
est une vari\'et\'e \'eventuellement singuli\`ere et avec ou sans
bord. 
\par
Soient $z\in J$ et $w\in \varphi^{-1}(z)$. Comme $\Lambda$ est
dilatante, la suite de $w_m:=\Lambda^{-m}(w)$ tend vers $0$. Par
cons\'equent la suite de $z_m:=\varphi(w_m)$ tend vers $b$. On a
$$g^m(z_m)=g^m(\varphi(w_m))=\varphi(\Lambda^m(w_m))=\varphi(w)=z.$$
D'o\`u $z_n\in g^{-n}(z)$. On en d\'eduit que $b$ appartient \`a
l'adh\'erence de $\bigcup_{m\geq 0} g^{-m}(z)$. Ceci est vrai
pour tout point p\'eriodique r\'epulsif $b$ de $f$
qui n'appartient \`a
la singularit\'e de $J$. L'ensemble de tels $b$ est dense dans
$J$. On constate que $\bigcup_{m\geq 0}g^{-m}(z)$ est dense dans
$J$.
\par
L'ensemble $\B_w$ est \'egal \`a $\varphi^{-1}(\bigcup_{m\geq
0}g^{-m}(z))$ o\`u $z:=\varphi(w)$.
Il est donc dense dans $\varphi^{-1}(J)=J^*$.
\end{preuve}
\begin{corollaire}
{\bf a.}  La valeur propre $\lambda_k$ de $\Lambda$ est
\'egale \`a $d_g$. La mesure $\mu^*$ est \'egale \`a
$c\Pi^*(m)$ o\`u
$c>0$ est une constante,
$\Pi: J^*\longrightarrow \C^{k-1}\times \R$ est d\'efinie par
$\Pi(z):=(z',\Re z_k)$ et $m$ est la mesure de Lebesgue sur
$\C^{k-1}\times \R$.
\par
{\bf b.} Il existe une application $\Lambda_f \in \G_4$ dont la
partie lin\'eaire est
de d\'eterminant $d^{(k+1)/2}$ telle que
$f\circ\varphi=\varphi\circ\Lambda_f$.
\end{corollaire}
\begin{preuve} {\bf a.} Posons $\tilde \mu^*:=\Pi^*(m)$.
D'apr\`es le  corollaire 3.2, on a 
  $\Lambda^*(\tilde \mu^*)=(\lambda_k)^k \mu^*$.
D'apr\`es le corollaire
  3.3, $\tilde\mu^*$ est invariante par $\A|_{J^*}$. D'apr\`es la
  proposition 3.4, il existe une mesure positive
$\tilde \mu$ \`a support
  dans $J$ telle que $\varphi^*(\tilde\mu)=\tilde\mu^*$. La relation
  $\Lambda^*(\tilde \mu^*)=\lambda_k^k \mu^*$ entra\^{\i}ne
  $g^*(\tilde \mu)=\lambda_k^k\tilde \mu$.
Mais on sait que $g$ d\'efinit un rev\^etement ramifi\'e de
degr\'e $(d_g)^k$
de $\C^k$ dans $\C^k$. Par
  cons\'equent la masse de  $g^*(\tilde \mu)$ est $(d_g)^k$ fois plus
  grande que celle de $\tilde \mu$. On obtient donc
  $(\lambda_k)^k=(d_g)^k$ et $\lambda_k=d_g$.
\par
On sait que $\mu$ est la seule mesure de probabilit\'e 
invariante de $g$ qui ne charge
pas les ensembles pluripolaires. La mesure $\tilde\mu$ est invariante
par $g$ et ne charge par les ensembles pluripolaires. Par
cons\'equent, il existe
une constante $c>0$ telle que $\mu=c\tilde\mu$. On en d\'eduit que
$\mu^*=c\Pi^*(m)$.
\par
{\bf b.} 
On choisit $w_1\in J^*$ et $w_2\in J^*$ tels que
$f(\varphi(w_1))=\varphi(w_2)$. On peut supposer que $w_1$, $w_2$
n'appartiennent pas \`a ${\cal C}_\varphi$ et $\varphi(w_1)$
n'appartient pas \`a ${\cal C}_f$. Fixons des petits voisinages
$W_1$, $U_1$, $W_2$, $U_2$ de $w_1$, $\varphi(w_1)$, $w_2$ et
$\varphi(w_2)$ tels que $\varphi$ r\'ealise un biholomorphisme entre
$W_1$ (resp $W_2$) et $U_1$ (resp. $U_2$) et $f$ r\'ealise un
biholomorphisme entre $U_1$ et $U_2$. Posons
$\Lambda_f:=(\varphi|_{W_2})^{-1}\circ (f|_{U_1})\circ
(\varphi|_{W_1})$. On a $f\circ\varphi=\varphi\circ \Lambda_f$.
Comme $f(J)=J$, on a $\Lambda(J^*\cap W_1)=J^*\cap
W_2$. D'apr\`es le lemme 3.1, $\Lambda_f\in \G_4$. Par prolongement
analytique $f\circ\varphi=\varphi\circ \Lambda_f$ sur $\C^k$. On
en d\'eduit que $g\circ\varphi=\varphi\circ \Lambda_f^n$ car
$g=f^n$. Ceci montre
que $\Lambda\circ \Lambda_f^{-n}$ pr\'eserve les fibres de
$\varphi$. Il est donc un \'el\'ement de $\A$ et ainsi
sa partie lin\'eaire est de d\'eterminant 1. 
Par cons\'equent, $\det
\Lambda_f^n=\det \Lambda=d_g^{(k+1)/2}=d^{n(k+1)/2}$. D'o\`u $\det
\Lambda_f=d^{(k+1)/2}$. 
\end{preuve}
\section{Fibration invariante}
\ahead
Observons que $\Lambda_f$ pr\'eserve la fibration $\{z'=\const\}$. On
montrera que l'image de cette fibration par $\varphi$ est une
fibration de droites complexes passant par un point. Cette fibration
est donc invariante par $f$ et ceci impliquera que $f$ est homog\`ene.
\par
Pour la suite, on note $L:=\P^k\setminus \C^k$ l'hyperplan \`a
l'infini. Pour tout $u\in \C^{k-1}$, on note 
$\D_u:=\{z\in\C^k, z'=u\}$, $\D_u^+:=\{z\in \D, \Im z_k\geq
\|u\|^2\}$, $\D_u^-:=\{z\in \D, \Im z_k\leq \|u\|^2\}$ et $l_u:=\{z\in
\D, \Im z_k = \|u\|^2\}$.
\par
Rappelons la version suivante du principe de
Phragm\'en-Lindel\"off qui sera
utilis\'ee dans ce paragraphe:
\begin{lemme}{\bf \cite{DinhSibony}} Soit $v\geq 0$
une fonction continue
sous-harmonique d\'efinie sur le demi-plan
$\D^-:=\{z\in \C|\ \Im z\leq 0\}$ et
nulle sur $\R$. 
Supposons qu'il
existe un nombre r\'eel $\lambda>1$ tel que $v(\lambda z)=\lambda
v(z)$ pour
tout $z\in \D^-$. Alors $v(z)=-c\Im z$ o\`u $c\geq 0$ est
une constante.
\end{lemme} 
\begin{proposition}
La fonction $G^*(z)$ est \'egale \`a
$\max(0,c'(\|z'\|^2-\Im z_k))$ o\`u $c'>0$ est une constante.
Par cons\'equent, les 
vari\'et\'es complexes $\Sigma\subset\{G^*\not=0\}$
sur lesquelles
$G^*$ est pluriharmonique sont des ouverts des droites complexes
$\D_u$.
\end{proposition}
\begin{preuve} D'apr\`es la proposition 2.6, si $\Im z_k\geq
\|z'\|^2$ on a $G^*(z',z_k)=0$. Il reste \`a traiter le cas o\`u
$\Im z_k\leq \|z'\|^2$. 
\par
D'apr\`es le lemme 4.1, il existe une constante
$c'\geq 0$ telle que $G^*(0,z_k)=-c'\Im z_k$ lorsque $\Im z_k\leq
0$. On va montrer que $G^*(z',z_k)=c'(\|z'\|^2-\Im z_k)$ pour
$\Im z_k\leq \|z'\|^2$. Comme $G^*$ est continue, il suffit de le
montrer pour une famille dense de $z'$.
Soient $\B_0$ d\'efini dans le corollaire 3.5
et $w=(u,v)\in \B_0$.
Alors il existe $s\geq 0$ et $\tau\in\A$ tels que
$w=\Lambda^{-s}\circ\tau(0)$. Supposons que $\Im z_k\leq
\|u\|^2$. Les applications $\tau$ et
$\sigma_{1/d_g}\circ\Lambda$ pr\'eservant les vari\'et\'es
$\{\Im z_k-\|z'\|^2=\const\}$,
il existe $x\in\R$ tel
que $\tau^{-1}\circ \Lambda^s(u,z_k)=(0,x+i(d_g)^s(\Im
z_k-\|z'\|^2))$.
D'autre part $G^*\circ\Lambda=d_gG^*$ et $G^*\circ\tau=G^*$. On
obtient
$$G^*(u,z_k)=(d_g)^{-s}G^*(0,x+i(d_g)^s(\Im z_k-\|z'\|^2))=-c'(\Im
z_k-\|z'\|^2).$$
Ceci est valable pour tout $w=(u,v)\in \B_0$. D'apr\`es le corollaire
3.5, $\B_0$ est dense dans $J^*$. D'o\`u $G^*(z',z_k)=
c'(\|z'\|^2-\Im z_k)$ pour tout
$\Im z_k\leq \|z'\|^2$. La fonction $G^*$ est positive ou nulle
et non identiquement nulle. D'o\`u $c'>0$.
\end{preuve}
\begin{lemme} Soit $u\in\C^{k-1}$. Supposons que la restriction de
  $\varphi$ sur $\D_u$ n'est pas injective. Alors 
le groupe $\A_u:=\A|_{\D_u}$ agit transitivement sur les fibres de
$\varphi|_{\D_u}$. De plus, il existe
  $\alpha\in \R^+$ tel que
$$\A_u=\{(u,z_k)\mapsto (u,z_k+n\alpha) \mbox{ avec } n\in
\Z\}.$$
\end{lemme}
\begin{preuve} Montrons d'abord que $\A_u$ contient un
  \'el\'ement diff\'erent de l'identit\'e. 
Soient $w_1, w_2\in \D_u$ v\'erifiant
$\varphi(w_1)=\varphi(w_2)$. D'apr\`es la proposition 3.4, il existe
$\tau\in \A$ tel que $\tau(w_1)=\tau(w_2)$ et $w_1\not=w_2$. Comme
$\tau\in \A\subset \G_3$,
$\tau$ pr\'eserve la fibration $\{z'=\const\}$. Par cons\'equent,
$\tau$ pr\'eserve la droite $\D_u$ qui contient $w_1$ et $w_2$. On
en d\'eduit que  $\A_u$ contient $\tau|_{\D_u}$ qui est
diff\'erent de l'identit\'e car $\tau(w_1)=w_2\not=w_1$.
\par
Ce raisonnement implique \'egalement que $\A_u$ agit
transitivement sur les fibres de $\varphi|_{\D_u}$ car $\A$
agit transitivement sur les fibres de $\varphi$.
\par
Comme $\A$ est un sous-groupe discret de $\G_3$, 
le groupe  $\A_u$ est un groupe discret
d'isom\'etries complexe de $\D_u$.
Le fait que $J^*$ est invariant par $\A$ implique que  $\A_u$
pr\'eserve la droite r\'eelle $l_u:=\D_u\cap J^*=\{z'=u \mbox{ et } \Im
z_k=\|u\|^2\}$. Soit $\tau\in \A_u$. Alors il existe $x\in
\R$ tel que $\tau(u,z_k)=(u,\pm z_k+x)$. Montrons que $\tau(u,z_k)
\not =(u,-z_k+x)$. Supposons que c'est le cas. Alors comme
$G^*\circ\tau=G^*$, la proposition 4.2 entra\^{\i}ne $c'=0$ ce qui
est impossible.
\par
Alors les \'el\'ements de  $\A_u$ sont du type 
$(u,z_k)\mapsto (u,z_k+x)$. Comme c'est un groupe discret, 
il existe
  $\alpha\in \R^+$ tel que
$$\A_u=\{(u,z_k)\mapsto (u,z_k+n\alpha) \mbox{ avec } n\in
\Z\}.$$
\end{preuve}
\par
On dit que ${\cal C}\subset \P^k$ est {\it une courbe holomorphe
immerg\'ee} s'il existe une surface de Riemann $S$ et une
application holomorphe non constante $\psi$ de $S$ dans $\P^k$
telles que ${\cal C}=\psi(S)$.
\begin{lemme} Soient $u\in\C^{k-1}$ et $\omega\in L$. Supposons
que l'ensemble
${\cal C}:=\varphi(\D_u)\cup \{\omega\}$
est une courbe holomorphe immerg\'ee
dans $\P^k$. Supposons aussi qu'en $\omega$ au moins une
composante de ${\cal C}$ 
coupe $L$ transversalement. Alors
il existe $\alpha\in \R^+$ et $\nu_1,\nu_2\in
  \C^k$ tels que $\varphi(u,z_k)=\exp(2\pi i z_k/\alpha)\nu_1
+\nu_2$. En particulier, $\varphi(\D_u)$ 
est contenue dans une droite complexe.
\end{lemme}
\begin{preuve} Montrons d'abord que $\varphi|_{\D_u}$ n'est pas
  injective. Supposons que cette application soit injective. Alors
  $\varphi(\D_u)$ est un $\C$ immerg\'e dans $\C^k$. Par
  hypoth\`ese,  $\varphi(\D_u)\cup \{\omega\}$ est une courbe
  holomorphe immerg\'ee dans $\P^k$. Cette courbe est donc un $\P^1$. 
  Par cons\'equent, $\varphi(u,z_k)$ tend vers
  $\omega$ quand $z_k\rightarrow \infty$. C'est contradiction car
  l'image de $\D_u^+$ par
  $\varphi$ est
  contenue dans le compact $\{G=0\}$ de $\C^k$.
\par
On note $\alpha$ le nombre r\'eel positif d\'efini dans le
lemme 4.3. Soit
$\Psi:\D_u\longrightarrow \C^*$
d\'efinie par $\Psi(u,z_k):=\exp(2i\pi z_k/\alpha)$.
Alors il existe une
application holomorphe injective
$\psi:\C^*\longrightarrow \C^k$ telle que
$\varphi=\psi\circ\Psi$. Posons $\tilde G:=G\circ\psi$. On a
$G^*=\tilde G\circ \Psi$ et d'apr\`es la proposition 4.2
$$\tilde G(\zeta)=\max\left(0,c'\|u\|^2+\frac{c'\alpha}{2\pi}
\log|\zeta|\right)$$
pour tout $\zeta\in\C^*$. On sait $\{G=0\}$ est un compact de
$\C^k$. Par cons\'equent, $\psi$ est born\'ee au
voisinage de $0$. Elle se prolonge holomorphiquement en une
application de $\C$ dans $\C^k$. Quand $z$ tend vers l'infini,
la fonction $G$ a une croissance logarithmique. Ceci implique que
$\psi(\zeta)$ a une croissance polynomiale
quand $\zeta\rightarrow \infty$. Par suite,
$\overline{\varphi(\D_u)}$ est
une courbe alg\'ebrique qui coupe $L$ en un seul point. Par
hypoth\`ese, ce point
doit \^etre  $\omega$ et l'intersection de la courbe alg\'ebrique
avec $L$ est transversale. Cette courbe est donc une droite
projective. 
\par
Comme $\psi$ est injective, elle
est une application polynomiale de degr\'e
1. Soient $\nu_1\in\C^k$ et $\nu_2\in \C^k$ tels que
$\psi(\zeta)=\nu_1\zeta +\nu_2$. Alors
$\varphi(0,z_k)=\exp(2i\pi z_k/\alpha)\nu_1+\nu_2$.
\end{preuve}
\begin{proposition} Il existe $\alpha\in \R^+$, $\nu_2\in\C^k$ et une
  application holomorphe $\nu_1:\C^{k-1} \longrightarrow \C^k$ tels
  que $\varphi(z',z_k)=\nu_1(z')\exp(2i\pi z_k/\alpha)+\nu_2$. 
\end{proposition}
\begin{preuve} La restriction de $f$ \`a $L$ est une application
holomorphe de degr\'e $d\geq 2$. Il existe donc un point
$\omega\in L$ qui est p\'eriodique r\'epulsif pour $f|_L$.
On peut supposer que
c'est un point fixe r\'epulsif de $g|_L$. Alors en $\omega$, $g$
poss\`ede $k-1$ valeurs propres de module sup\'erieur \`a 1 et
une valeur propre nulle. Soit
$$V:=\{z\in \P^k, \lim_{s\rightarrow\infty} g^s(z)=\omega\}$$
la vari\'et\'e stable de $g$ en
$\omega$ \cite[p.27]{Ruelle}. C'est une courbe holomorphe immerg\'ee
dans $\P^k$ qui coupe $L$ transversalement en $\omega$. Montrons que 
$G$ est harmonique sur $V$. Notons $[w_1:\cdots:w_k]$ les
coordonn\'ees homog\`enes de $w$ dans $L$.
Sans nuire \`a la g\'en\'eralit\'e,
on suppose que $w_k=1$.
Soit $z\in V$. Notons $g_{s,m}(z)$ la $m$-i\`eme
coordonn\'ee de $g^s(z)$. Alors $\lim_{s\rightarrow\infty}
g^s(z)=w$ et $\lim_{s\rightarrow\infty}g_{s,m}(z)/g_{s,k}(z)=w_m$
pour tout $z\in V$ et tout $m$. On a
\begin{eqnarray*}
G(z) & = &\lim_{s\rightarrow\infty} \frac{1}{(d_g)^s}\log
\|g^s(z)\| \\
& = &
\lim_{s\rightarrow\infty} \frac{1}{(d_g)^s}\log |g_{s,k}(z)| +
\frac{1}{2(d_g)^s}\log\sum_{m=1}^k
\left |\frac{g_{s,m}(z)}{g_{s,k}(z)}
\right |^2  
\end{eqnarray*} 
Le second terme de la derni\`ere expression tend vers $0$. Le
premier terme tend vers une fonction harmonique car il est
harmonique. Ceci montre que $G$ est harmonique sur $V$
\par
D'apr\`es la proposition
4.2, il existe $\xi\in\C^{k-1}$ tel que $V\subset\varphi(\D_\xi)$.
D'apr\`es le lemme 4.4, $\varphi(\D_\xi)$ est une droite complexe.
\par
Soit $\Sigma$ l'ensemble des $u\in\C^{k-1}$ tels que $\D_u$
v\'erifie
$\Lambda^m(\D_u)=\tau(\D_\xi)$ pour un $m\geq 0$ et un
$\tau\in\A$. D'apr\`es le corollaire 3.5, cet ensemble est
dense dans $\C^{k-1}$. Pour un tel $u$, 
$g^m(\varphi(\D_u))$ est \'egal \`a
$\varphi(\D_\xi)$ qui est contenu dans une droite complexe.
Par cons\'equent,
$\overline{\varphi(\D_u)}$ est une courbe alg\'ebrique qui
coupe $L$ transversalement. D'apr\`es le lemme 4.4, c'est une
droite projective. Alors pour tout $u\in \Sigma$, il existe
$\alpha(u)\in\R^+$, $\nu_1(u)\in\C^k$ et $\nu_2(u)\in\C^k$ tels que
$\varphi(u,z_k)=\exp(2i\pi z_k/\alpha(u))\nu_1(u)+\nu_2(u)$.
\par
L'application $\varphi$ \'etant holomorphe,
en prenant $z_k=0,1,-1$, on montre facilement
que $\exp(1/\alpha(u))$, $\nu_1(u)$ et
$\nu_2(u)$ d\'ependent 
holomorphiquement de $u$. Comme $\alpha(u)$ prend des valeurs
r\'eelles dans une famille dense de $u$, elle est constante. On
obtient des formules explicites de $\nu_1(u)$ et $\nu_2(u)$ en
fonction de $\varphi$ qui, par
continuit\'e, impliquent que $\varphi(u,z_k)=\nu_1(u)\exp(2i\pi
z_k/\alpha) +\nu_2(u)$ pour tout $u$.
Comme $G^*$ est nulle sur $\{\Im
z_k\geq \|z'\|^2\}$, $G$ doit s'annuller au point
$\nu_2(u)=\lim_{\Im z_k \rightarrow +\infty}\varphi(u,z_k)$. Par
cons\'equent, 
$\nu_2(u)$ est une application holomorphe de $\C^{k-1}$
\`a l'image dans le
compact $\{G=0\}$ de $\C^k$. D'apr\`es le th\'eor\`eme de
Liouville, elle doit \^etre constante.
\end{preuve}
{\it Preuve du th\'eor\`eme 1.3---} Quitte \`a effectuer un changement
lin\'eaire de coordonn\'ees de $\C^k$,
on peut supposer $\nu_2=0$. L'application $\Lambda_f$
pr\'eserve la fibration $\{z'=\const\}$. D'apr\`es la proposition
4.5,
$f$ pr\'eserve la fibration des droites passant par $0$. Elle est donc
une application homog\`ene.
\par
Notons $\P^{k-1}$ l'ensemble des droites complexes passant par $0$
que l'on peut identifier avec l'hyperplan \`a l'infini de
$\C^k$. L'application $f$
induit un endomorphisme holomorphe $\widehat f$ de
$\P^{k-1}$. Il reste \`a prouver que $\widehat f$ est de Latt\`es. On
d\'efinit l'application $\widehat \varphi$ de
$\C^{k-1}$ dans $\P^{k-1}$
par $z'\mapsto [\nu_1(z')]$. Notons $\widehat\A$ la restriction des
actions de $\A$ sur la famille de droites $\{\D_u\}_{u\in\C^{k-1}}$ et
$\widehat\Lambda_f$ la restriction de l'action de
$\Lambda_f$ sur cette famille de droites.
On a $\widehat f\circ\widehat
\varphi=\widehat\varphi\circ\widehat \Lambda_f$ car
$f\circ\varphi=\varphi\circ\Lambda_f$. D'apr\`es la proposition 3.4,
$\widehat \A$ est un groupe
d'isom\'etries co-compact qui agit transitivement sur les fibres de
$\widehat\varphi$. 
Par d\'efinition, $\widehat
f$ est de Latt\`es. On en d\'eduit que la restriction de $f$ \`a
$L$ est aussi une application de Latt\`es.
\par
\hfill $\square$
\\
{\it Preuve du corollaire 1.5---} Soit $F$ un relev\'e polynomial de
$f$. Soient
$$G_f(z):=\lim_{s\rightarrow \infty}\frac{1}{d^s}\log \|F^s(z)\|$$
la fonction de Green de $f$ et $G^+(z):=\max(0,G_f(z))$ la
fonction de Green de $F$. Notons $\mu_F$ la mesure
d'\'equilibre de $F$ et
$\pi:\C^{k+1}\setminus\{0\}\longrightarrow \P^k$ est la
projection canonique.
\par
Comme le courant $T$ de $f$ est de classe
${\cal C}^{2+\epsilon}$ et
est strictement positif dans un ouvert $U$, la fonction
$G_f$ est de classe ${\cal C}^{2+\epsilon}$
dans $\pi^{-1}(U)$. En particulier, $\{G_f=0\}\cap \pi^{-1}(U)$
est une
hypersurface ${\cal C}^{2+\epsilon}$ strictement pseudoconvexe.
Comme $G_f(\lambda z)=\log|\lambda|+G_f(z)$,
la fonction $G_f$ est harmonique sur toute droite complexe passant
par $0$. On en d\'eduit que dans $\pi^{-1}(U)$ le support de
$\mu_F=(\ddc G^+)^{k+1}$ est \'egale \`a $\{G_f=0\}$.
D'apr\`es le
th\'eor\`eme 1.3 appliqu\'e \`a $F$,
la restriction de $F$ \`a l'infini est une
application de
Latt\`es. Ceci implique que $f$ est \'egalement
une application de Latt\`es.
Tien-Cuong Dinh\\
Math\'ematique - B\^atiment 425 \\
UMR 8628, Universit\'e Paris-Sud \\
91405 ORSAY Cedex (France) \\
TienCuong.Dinh@math.u-psud.fr
\end{document}